\documentclass[MSNbibl,nameyear,seceqn,dvips]{arxstspdf}
\usepackage{mathbh}
\usepackage{graphicx}
\usepackage{flushend}
\usepackage{stfloats}

\volume{27}
\issue{4}
\pubyear{2012}
\firstpage{558}
\lastpage{575}
\doi{10.1214/12-STS393} 

\makeatletter
\def\eqref#1{(\ref{#1})}
\newcommand{\pP}{\mathsf{P}}
\newcommand{\qQ}{\mathsf{Q}}
\newcommand{\D}{R}
\newcommand{\cF}{\mathcal{F}}
\newcommand{\cG}{\mathcal{G}}
\newcommand{\cH}{\mathcal{H}}
\newcommand{\cK}{\mathcal{K}}
\newcommand{\cL}{\mathcal{L}}
\newcommand{\cM}{\mathcal{M}}
\newcommand{\cN}{\mathcal{N}}
\newcommand{\cP}{\mathcal{P}}
\newcommand{\cX}{\mathcal{X}}
\newcommand{\cZ}{\mathcal{Z}}
\newcommand{\bI}{\mathbf{I}}
\newcommand{\bX}{\mathbf{X}}
\newcommand{\bY}{\mathbf{Y}}
\newcommand{\R}{\mathbb{R}}
\newcommand{\E}{\mathbb{E}}
\newcommand{\1}{\mathbh{1}}
\newcommand{\X}{\mathbb{X}}
\newcommand{\qqq}{\mathsf{q}}
\newcommand{\pp}{\mathsf{p}}
\newcommand{\ff}{\mathsf{f}}
\newcommand{\Tr}{\operatorname{Tr}}
\newcommand{\argmin}{\mathop{\operatorname{argmin}}}
\newtheorem{TH1}{Theorem}[section]
\newtheorem{cor}{Corollary}[section]
\newtheorem{lem}{Lemma}[section]
\newproclaim{Fusedsparsity}{Fused sparsity}
\newcommand{\bookbox}[1]{
\par\medskip\noindent
\framebox[\columnwidth]{
\begin{minipage}{235pt}
{#1}
\end{minipage} } \par\medskip }
\makeatother

\begin{document}
\begin{frontmatter}

\title{Sparse Estimation by Exponential Weighting}
\runtitle{Sparsity by exponential weighting}

\begin{aug}
\author[a]{\fnms{Philippe} \snm{Rigollet}\corref{}\ead[label=e1]{rigollet@princeton.edu}}
\and
\author[b]{\fnms{Alexandre B.} \snm{Tsybakov}\ead[label=e2]{alexandre.tsybakov@ensae.fr}}
\runauthor{P. Rigollet and A. B. Tsybakov}

\affiliation{Princeton University and CREST-ENSAE}

\address[a]{Philippe Rigollet is Assistant Professor, Department of Operations
Research and Financial Engineering, Princeton University, Princeton,
New Jersey 08544, USA \printead{e1}.}
\address[b]{Alexandre B. Tsybakov is Professor and Head, Laboratoire de Statistique,
CREST-ENSAE,
3, av. Pierre Larousse, F-92240 Malakoff Cedex, France
\printead{e2}.}

\end{aug}

%
\begin{abstract}
Consider a regression model with fixed design and Gaussian noise where
the regression function can potentially be well approximated by a
function that admits a sparse representation in a given dictionary.
This paper resorts to exponential weights to exploit this underlying
sparsity by implementing the principle of \textit{sparsity pattern
aggregation}. This model selection take on sparse estimation allows us
to derive sparsity oracle inequalities in several popular frameworks,
including ordinary sparsity, fused
sparsity and group sparsity. One striking aspect of these theoretical
results is that they hold under \textit{no condition in the dictionary}.
Moreover, we describe an efficient
implementation of the sparsity pattern aggregation principle that
compares favorably to state-of-the-art procedures on some basic
numerical examples.
\end{abstract}

%

\end{frontmatter}

\section{Introduction}
\label{sec1}

Since the 1990s, the idea of exponential weighting has been
successfully used in a variety of statistical problems. 
In this paper, we review several properties of estimators based on
exponential weighting with a particular emphasis on how they can be
used to construct optimal and computationally efficient procedures
for high-dimensional regression under the sparsity scenario.

Most of the work on exponential weighting deals with a regression
learning problem. Some of the results can be extended to other
statistical models such as density estimation or classification; cf.
Section~\ref{SEC:related}. For the sake of brevity and to make the
presentation more transparent, we focus here on the following
framework considered in \citet{RigTsy11}. Let $\cZ =\{(x_1, Y_1),
\ldots, (x_n, Y_n)\}$ be a collection of independent random pairs
such that $(x_i,Y_i) \in \cX \times \R$, where $\cX$ is an arbitrary
set. Assume the regression model
\begin{equation}
\label{EQ:model}
Y_i=\eta(x_i)+\xi_i, \quad i=1, \ldots, n ,
\end{equation}
where $\eta\dvtx\cX \to \R$ is the unknown regression function, and the
errors $\xi_i$ are independent Gaussian $\cN(0,\sigma^2)$. The
covariates are deterministic elements $x_1, \ldots, x_n$ of~$\cX$.
For any function $f\dvtx\cX \to \R$, we define a seminorm $\|\cdot\|$
by\footnote{Without loss
of generality, in what follows we will associate all the functions
with vectors in $\R^n$ since only the values of functions at points
$x_1,\dots,x_n$ will appear in the risk. So, $\|\cdot\|$ will be
indeed a norm and, with no ambiguity, we will use other related
notation such as $\|\bY - f\|$ where $\bY$ is a vector in $\R^n$
with components $Y_1, \ldots, Y_n$.}
\[
\|f\|^2 =\frac{1}{n}\sum_{i=1}^n f^2(x_i) .
\]
%
%
%

We adopt the following learning setup. 
Let $\cH=\{f_1, \ldots, f_M\}$, be a dictionary of $M\ge 1$ given
functions. For example, $f_j$ can be
some basis functions or some preliminary estimators of $f$
constructed from another sample that we consider as frozen; see
Section~\ref{SEC:estim} for more details. Our goal is to approximate the
regression function $\eta$ by a linear combination ${\sf
f}_\theta(x) = \sum_{j=1}^M \theta_j f_j(x)$ with weights $\theta =
(\theta_1,\dots,\theta_M)$, where possibly $M\gg n$. The performance
of a given estimator $\hat f$ of a function $\eta$ is measured in
terms of its averaged squared error
\[
R(\hat f)=\|\hat f -\eta\|^2:=\frac{1}{n}\sum_{i=1}^n [\hat
f(x_i) - \eta(x_i) ]^2 .
\]
Let $\Theta$ be a given subset of $\R^M$. In the aggregation
problem, we would ideally wish to find an \textit{aggregated estimator}
$\hat f$ whose risk $R(\hat f)$ is as close as possible in a
probabilistic sense to the minimum risk $\inf_{\theta \in \Theta}
R(\ff_\theta)$. Namely, one can construct
estimators $\hat f$ satisfying the following property:
\begin{equation}
\label{EQ:standardOI}
\E R(\hat f) \le C\inf_{\theta \in \Theta}
R(\ff_\theta) + \delta_{n,M}(\Theta) ,
\end{equation}
where $\delta_{n,M}(\Theta)$ is a small remainder term
characterizing the performance of the given aggregate $\hat f$ and
the complexity of the set $\Theta$, $C\ge1$ is a constant, and $\E$
denotes the expectation. Bounds of the form \eqref{EQ:standardOI}
are called \textit{oracle inequalities}. In some cases, even more
general results are available. They have the form
\begin{equation}
\label{EQ:balOI}
\E R(\hat f) \le C\inf_{\theta \in \Theta'} \{
R(\ff_\theta) +\Delta_{n,M}(\theta) \} ,
\end{equation}
where $\Delta_{n,M}$ is a remainder term that characterizes the
performance of the given aggregate $\hat f$ and the complexity of
the parameter $\theta \in \Theta'\subseteq \R^M$ (often \mbox{$\Theta'=
\R^M$}). To distinguish from \eqref{EQ:standardOI}, we will call
bounds of the form \eqref{EQ:balOI} the \textit{balanced oracle
inequalities}. If $\Theta \subseteq \Theta'$, then
\eqref{EQ:standardOI} is a direct consequence of \eqref{EQ:balOI}
with $\delta_{n,M}(\Theta)=C\sup_{\theta \in \Theta}
\Delta_{n,M}(\theta)$.

In this paper, we mainly focus on
the case where the complexity of a vector $\theta$ is measured as
the number of its nonzero coefficients $|\theta|_0$. In this case,
inequalities of the form~\eqref{EQ:balOI} are sometimes called
\textit{sparsity oracle inequalities}. Other measures of complexity,
also related to \textit{sparsity} are considered in
Section~\ref{sub:priors}. As indicated by the notation and
illustrated below, the remainder term $\Delta_{n,M}(\theta)$ depends
explicitly on the size $M$ of the dictionary and the sample size
$n$. It reflects the interplay between these two fundamental
parameters and also the complexity of~$\theta$.\vadjust{\goodbreak}


When the linear model is misspecified, that is, where there is no
$\theta \in \Theta$ such that $\eta=\ff_\theta$ on the\break set
$\{x_1,\dots,x_n\}$, the minimum risk satisfies\linebreak[4] $\inf_{\theta \in
\Theta} R(\ff_\theta)>0$ leading to a systematic bias term. Since
this term is unavoidable, we wish to make its contribution as small
as possible, and it is therefore important to obtain a leading
constant $C=1$. Many oracle inequalities with leading constant $C>1$
can be found in the literature for related problems. However, in
most of the papers, the set $\Theta=\Theta_n$ depends on the sample
size $n$ in such a way that $\inf_{\theta \in \Theta_n}R(\ff_\theta)$
tends to
$0$ as $n$ goes to infinity, under additional regularity
assumptions. In this paper,
we are interested in the case where $\Theta$ is fixed. 
For this reason, we consider here only oracle inequalities with
leading constant $C=1$ (called \textit{sharp oracle inequalities}).
Because they hold for finite $M$ and $n$, these are truly finite
sample results.

One salient feature of the oracle approach as opposed to standard
statistical reasoning, is that it does not rely on an underlying
model. Indeed, the goal is not to estimate the parameters of an
underlying ``true'' model but rather to construct an estimator that
mimics, in terms of an appropriate oracle inequality, the
performance of the best model in a given class, whether this model
is true or not. From a statistical viewpoint, this difference is
significant since performance cannot be evaluated in terms of
parameters. Indeed, there is no true parameter. However, we can
still compare the risk of the estimator with the optimum value.
Oracle inequalities offer a tool for such a comparison.

A particular choice of $\Theta$ corresponds to the problem of
\textit{model selection aggregation}. Let $\Theta=\Theta^{\rm MC}$ be the set
of $M$ canonical basis vectors of $\R^M$. Then the set of linear
combinations $\{\ff_\theta, \theta\in \Theta^{\rm MC}\}$ coincides
with the initial dictionary of functions $\cH=\{f_1, \ldots, f_M\}$,
so that the goal of model selection is to mimic the best function in
the dictionary in the sense of the risk measure $R(\cdot)$. This can
be done in different ways, leading to different rates
$\delta_{n,M}(\Theta^{\rm MC})$; however one is mostly interested in
the methods that attain the rate $\delta_{n,M}^*(\Theta^{\rm
MC})\asymp(\log M)/n$ which is known to be minimax optimal
(see \cite{Tsy03}; \cite{BunTsyWeg07c}; Rigollet\break (\citeyear{Rig12})). The first sharp oracle
inequalities with this rate for a setting different from the one
considered here were obtained by \citet{Cat99} (see also \cite{Cat04}), who used the progressive mixture method based on
exponential weighting. Other methods of model selection for
aggregation consist in selecting a function in the dictionary by
minimizing a (penalized) empirical risk (see, e.g., \cite{Nem00}; \cite{Weg03}; \cite{Tsy03}; \cite{Lec11}). One of the major novelties
offered by exponential weighting is to \textit{combine} (average) the
functions in the dictionary using a convex combination, and not
simply to \textit{select} one of them. From the theoretical point of
view, selection of one of the functions has a fundamental drawback
since it does not attain the optimal rate $(\log M)/n$; cf.
Section~\ref{SEC:exp}. 

The rest of the paper is organized as follows. In the next section,
we discuss some connections between the exponential weighting
schemes and penalized empirical risk minimization.
In Section~\ref{SEC:OI}, we present the first oracle inequalities
that demonstrate how exponential weighting can be used to
efficiently combine functions in a dictionary. The results of
Section~\ref{SEC:OI} are then extended to the case where one wishes
to combine not deterministic functions, but
estimators. Oracle inequalities for this problem are discussed in
Section~\ref{SEC:estim}. They are based on the work of
\citet{LeuBar06} and \citet{DalSal11}. 
Section~\ref{SEC:sparse} shows how these results can be adapted to
deal with sparsity. We introduce the principle of \textit{sparsity
pattern aggregation}, and we derive sparsity oracle inequalities in
several popular frameworks including ordinary sparsity, fused
sparsity and group sparsity. Finally, we describe an efficient
implementation of the sparsity pattern aggregation principle and
compare its performance to state-of-the-art procedures on some basic
numerical examples.

\section{Exponential Weighting and
Penalized Risk Minimization}

\label{SEC:exp}

\subsection{Suboptimality of Selectors}

A natural candidate to solve the problem of model selection introduced in the previous section is an empirical risk minimizer. Define the empirical risk by
\[
\hat R_n(f)=\frac{1}{n}\sum_{i=1}^n  [Y_i -
f(x_i) ]^2=\|\bY-f\|^2
\]
and the empirical risk minimizer by
\begin{equation}
\label{EQ:defERM}
\hat f^{\textsc{erm}}
=\argmin_{f \in \cH}\hat R_n(f),
\end{equation}
where ties are broken arbitrarily. However, while this procedure
satisfies an exact oracle inequality, it fails to exhibit the
optimal rate of order $\delta_{n,M}^*(\Theta^{\rm MC})\asymp(\log
M)/n$.\vadjust{\goodbreak} The following result shows that this defect is intrinsic not
only to empirical risk minimization but also to any method that
selects only one function in the dictionary~$\cH$. This includes
methods of model selection by penalized empirical risk minimization.
We call estimators $\hat S_n$, taking values in $\cH$ the
\textit{selectors}.

\begin{TH1}
\label{TH:ERM}
Assume that $\|f_j\|\le 1$ for any \mbox{$f_j \in \cH$}. Any empirical risk minimizer
$\hat f^\textsc{erm}$ defined in~\eqref{EQ:defERM} satisfies the following oracle inequality:
\begin{equation}
\label{EQ:TH:ERM:UB}
\qquad\E R(\hat f^\textsc{erm}) \le \min_{1\le j\le
M} R(f_j) + 4\sigma \sqrt{\frac{2\log M}{n}}.
\end{equation}
Moreover, assume that
\begin{equation}
\label{EQ:C0}
(\sigma\vee 1)\sqrt{(\log M)/n} \le C_0
\end{equation}
for $0<C_0<1$ small enough. Then, there exists a dictionary
$\cH=\{f_1,\dots,f_M\}$ with  $\|f_j\|\le 1$, $j=1,\dots,M,$ such
that the following holds. For any selector $\hat S_n$, and in
particular, for any selector based on penalized empirical risk
minimization, there exists a regression function $\eta$ such that
$\|\eta\| \le 1$ and
\begin{equation}
\label{EQ:TH:ERM:LB} \E R(\hat S_n) \ge \min_{1\le j\le M} R(f_j) +
C_*\sigma\sqrt{\frac{\log M}{n}}
\end{equation}
for some positive constant $C_*$.
\end{TH1}

\begin{pf}
See the \hyperref[app]{Appendix}.
\end{pf}

It follows from the lower bound \eqref{EQ:TH:ERM:LB} that
\textit{selecting} one of the functions in a finite dictionary $\cH$
to solve the problem of model selection is suboptimal in the sense
that it exhibits a too large remainder term, of the order
$\sqrt{(\log M)/n}$. It turns out that we can do better if we take a
\textit{mixture}, that is, a convex combination of the functions
in~$\cH$. We will see in Section~\ref{SEC:OI} [cf.~\eqref{modsel_ora1}] that under a particular choice of weights in
this convex combination, namely the \textit{exponential weights}, one
can achieve oracle inequalities with much better rate $(\log M)/n$.
This rate is known to be optimal in a minimax sense in several
regression setups, including the present one
(see \cite{Tsy03}; \cite{BunTsyWeg07c}; \cite{Rig12}).

\subsection{Exponential Weighting as a Penalized Procedure}

Penalized empirical risk minimization for model selection has
received a lot of attention in the literature, and many choices for
the penalty can be considered (see, e.g., \cite{BirMas01}; \cite{BarBouLug02}; \cite{Weg03};
\cite{LugWeg04}; \cite{BunTsyWeg07c}) to obtain oracle\vadjust{\goodbreak}
inequalities with the optimal or near optimal remainder term.
However, all these inequalities exhibit a constant $C>1$ in front of
the leading term. This is not surprising as we have proved in the
previous section that it is impossible for selectors to satisfy
oracle inequalities like~\eqref{EQ:standardOI} that are both sharp
(i.e., with $C=1$) and have the optimal remainder term. To overcome
this limitation of selectors, we look for estimators obtained as
convex combinations of the functions in the dictionary.

The coefficients of convex combinations belong to the flat simplex
\[
\Lambda^M:= \Biggl\{ \lambda \in \R^M  \dvtx  \lambda_j\ge 0,
\sum_{j=1}^M \lambda_j=1 \Biggr\}.
\]
Let us now examine a few ways to obtain potentially good convex
combinations. One candidate is a solution of the following penalized
empirical risk minimization problem:
\[
\min_{\lambda \in \Lambda^M}  \{ \hat R_n(\ff_\lambda)
+\mathrm{pen}(\lambda) \},
\]
where $\mathrm{pen}(\cdot)\ge 0$ is a penalty function. This choice
looks quite natural since it provides a proxy of the right-hand side
of the oracle inequality \eqref{EQ:balOI} where the unknown risk
$R(\cdot)$ is replaced by its empirical counterpart $\hat
R_n(\cdot)$. The minimum is taken over the simplex $\Lambda^M$
because we are looking for a convex combination. Clearly, the
penalty $\mathrm{pen}(\cdot)$ should be carefully chosen and ideally
should match the best remainder term $\Delta_{n,M}(\cdot)$.  Yet,
this problem may be difficult to solve as it involves a minimization
over $\Lambda^M$. Instead, we propose to solve a simpler problem.
Consider the following linear upper bound on the empirical risk:
\[
\sum_{j=1}^M \lambda_j \hat R_n(f_j)\ge \hat R_n(\ff_\lambda)\quad
  \forall  \lambda \in \Lambda^M
\]
and solve the following optimization problem:
\begin{equation}\label{opt_prob}
\min_{\lambda \in \Lambda^M}  \Biggl\{\sum_{j=1}^M
\lambda_j \hat R_n(f_j) +\mathrm{pen}(\lambda) \Biggr\}.
\end{equation}
Note that if $\mathrm{pen}(\lambda)\equiv 0$, the solution $\hat
\lambda$ of~\eqref{opt_prob} is simply the empirical risk minimizer
over the vertices of the simplex so that $\ff_{\hat \lambda}=\hat
f^{\textsc{erm}}$. In general, depending on the penalty function,
this problem may be more or less difficult to solve. It turns out
that the Kullback--Leibler penalty leads to a particularly simple
solution and allows us to approximate the best remainder term
$\Delta_{n,M}(\cdot)$ thus adding great flexibility to the resulting
estimator.\vadjust{\goodbreak}

Observe that vectors in $\Lambda^M$ can be associated to probability
measures on $\{1, \ldots, M\}$. Let $\lambda=(\lambda_1,\allowbreak  \ldots,
\lambda_M)$ and $\pi=(\pi_1,\ldots, \pi_M)$ be two probability
measures on $\{1, \ldots, M\}$,  and define the Kullback--Leibler
divergence between $\lambda$ and $\pi$ by
\[
\cK(\lambda, \pi)=\sum_{j=1}^M \lambda_j \log \biggl(
\frac{\lambda_j}{\pi_j} \biggr)\ge 0.
\]
Here and in the sequel, we adopt the convention that $0\log 0=0$,
$0\log (a/0)=0$, and $\log (a/0)=\infty$,  for any $a>0$.

Exponential weights can be obtained as the solution of the following
minimization problem. Fix $\beta>0$, a \textit{prior} $\pi \in
\Lambda^M$, and define the vector $\hat \lambda^\pi$ by
\begin{equation}\label{opt_prob1}
\qquad\hat \lambda^\pi=\argmin_{\lambda \in \Lambda^M} \Biggl\{\sum_{j=1}^M
\lambda_j \hat R_n(f_j) + \frac{\beta}{n}\cK(\lambda,
\pi) \Biggr\}.
\end{equation}
This constrained convex optimization problem has a unique solution
that can be expressed explicitly. Indeed, it follows from the
Karush--Kuhn--Tucker (KKT) conditions that the components $\hat
\lambda_j^\pi$ of $\hat\lambda^\pi$ satisfy
\begin{eqnarray}
\label{EQ:KKT}
n\hat R_n(f_j) + \beta\log \biggl( \frac{\hat
\lambda^\pi_j}{\pi_j} \biggr) + \mu  -\delta_j=0,\\
\eqntext{j=1, \ldots,M,}
\end{eqnarray}
where $\mu, \delta_1, \ldots, \delta_M\ge0$ are Lagrange
multipliers,
and
\[
\hat \lambda^\pi_j \ge 0, \quad  \delta_j \hat \lambda^\pi_j=0,\quad
  \sum_{j=1}^M \hat \lambda^\pi_j =1.
\]
Equation~\eqref{EQ:KKT} together with the above constraints lead to 
the following closed form solution:
\begin{eqnarray}
\label{EQ:defexp}
 \hat \lambda^\pi_j= \frac{\exp(-n\hat
R_n(f_j)/\beta)\pi_j} {\sum_{k=1}^M\exp(-n\hat
R_n(f_k)/\beta)\pi_k},\\
\eqntext{j=1, \ldots, M,}
\end{eqnarray}
called the \textit{exponential weights}.
We see that one immediate effect of penalizing by the
Kullback--Leibler divergence is that the solution
of~\eqref{opt_prob1} is not a selector. As a result, it achieves the
desired effect of averaging as opposed to selecting.



%
%

%


\section{Oracle Inequalities}
\label{SEC:OI}

An \textit{aggregate} is an estimator
defined as a weighted average of the functions in the
dictionary~$\cH$ with some\vadjust{\goodbreak} data-dependent weights. We focus on the
aggregate with exponential weights,
\[
\hat f^\pi = \sum_{j=1}^M \hat \lambda^\pi_j f_j,
\]
where $\hat \lambda^\pi_j$ is given in~\eqref{EQ:defexp}. This
estimator satisfies the following oracle inequality.

\begin{TH1}
\label{TH:PAC1}
The aggregate $\hat f^\pi$ with $\beta\ge 4\sigma^2$ satisfies the following balanced oracle inequality
\begin{equation}
\label{EQ:PAC1}
\hspace*{15pt}\E R(\hat f^\pi) \le \min_{\lambda \in
\Lambda^M} \Biggl\{ \sum_{j=1}^M \lambda_j R(f_j) + \frac{\beta}{n}
\cK(\lambda, \pi) \Biggr\}.\hspace*{-15pt}
\end{equation}
\end{TH1}

Comparing with (\ref{opt_prob1}) we see that $\hat \lambda^\pi$ is the minimizer of the unbiased estimator of the right-hand side of (\ref{EQ:PAC1}).
The proof of Theorem~\ref{TH:PAC1} can be found in the papers of Dalalyan and Tsybakov (\citeyear{DalTsy07,DalTsy08})
containing more general results. In particular, they apply to
non-Gaussian distributions of errors $\xi_i$ and to exponential
weights with a general (not necessarily discrete) probability
distribution $\pi$ on $\R^M$. Dalalyan and Tsybakov (\citeyear{DalTsy07,DalTsy08}) show that the
corresponding exponentially weighted aggregate $\hat f^\pi_*$ satisfies
the following bound:
\begin{equation}
\label{EQ:PAC2}
\hspace*{12pt}\E R(\hat f^\pi_*) \le \inf_{p} \biggl\{  \int
R(\ff_\theta) p({\mathrm{d}}\theta) + \frac{\beta}{n} \cK(p,
\pi) \biggr\},\hspace*{-12pt}
\end{equation}
where the infimum is taken over all probability distributions $p$ on
$\R^M$, and $\cK(p, \pi)$ denotes the Kullback--Leibler divergence
between the general probability measures $p$ and $\pi$. Bound
\eqref{EQ:PAC1} follows immediately from~\eqref{EQ:PAC2} by taking
$p$ and $\pi$ as discrete distributions.


A useful consequence of~\eqref{EQ:PAC1} can be obtained by
restricting the minimum on the right-hand side to the vertices of
the simplex~$\Lambda^M$. These vertices are precisely the vectors
$e^{(1)}, \ldots, e^{(M)}$ that form the canonical basis of $\R^M$
so that
\[
\sum_{j=1}^M e^{(k)}_j R(f_j) =R(f_k),
\]
where $e^{(k)}_j=\delta_{jk}$ is the $j$th coordinate of $e^{(k)}$,
with $\delta_{jk}$ denoting the Kronecker delta. It yields
\begin{equation}\label{modsel_ora}
\quad\E R(\hat f^\pi) \le \min_{1\le j \le M} \biggl\{  R(f_j) +
\frac{\beta}{n} \log (\pi_j^{-1}) \biggr\}.
\end{equation}
Taking $\pi$ to be the uniform distribution on $\{1, \ldots, M\}$
leads to the following oracle inequality:
\begin{equation}\label{modsel_ora1}
\E R(\hat f^\pi) \le \min_{1\le j \le M}  R(f_j) + \frac{\beta\log
M}{n} ,
\end{equation}
that exhibits a remainder term of the optimal order $(\log M)/n$.\vadjust{\goodbreak}


The role of the distribution $\pi$ is to put a prior weight on the
functions in the dictionary. When there is no preference, the
uniform prior is a common choice. However, 
we will see in Section~\ref{SEC:sparse} that choosing nonuniform
weights depending on suitable \textit{sparsity} characteristics can be
very useful. Moreover, this methodology can be extended to many
cases where one wishes to learn with a prior. It is worth mentioning
that while the terminology is reminiscent of a Bayesian setup, this
paper deals only with a frequentist setting (the risk is not
averaged over the prior).

\section{Aggregation of Estimators}
\label{SEC:estim}

\subsection{From Aggregation of Functions to Aggregation of Estimators}

Akin to the setting of the previous section, exponential weights
were originally introduced to aggregate deterministic functions
$f_j$ from a dictionary. These functions can be chosen in
essentially two ways. Either they have good approximation properties
such as an (over-complete) basis of functions or they are
constructed as preliminary estimators using a hold-out sample. The
latter case corresponds to the problem of \textit{aggregation of
estimators} originally described in \citet{Nem00}. The idea put
forward by \citet{Nem00} is to obtain two independent samples from
the initial one by randomization; estimators are constructed from the first sample while the second is used to perform
aggregation. To carry out the analysis of the aggregation step, it
is enough to work conditionally on the first sample so that the
problem reduces to aggregation of deterministic functions.
A limitation is that Nemirovski's randomization only applies to
Gaussian model with known variance. Nevertheless, this idea of two-step procedures
carries over to
models with i.i.d. observations where one can do direct sample
splitting (see, e.g., \cite{Yan04}; \cite{RigTsy07}; \cite{Lec07b}). Thus, in
many cases aggregation of estimators can be achieved by reduction to
aggregation of functions.

Along with this approach, one can aggregate estimators using the
same observations for both estimation and aggregation. While for
general estimators this would clearly result in overfitting,
the idea proved to be successful for
certain types of estimators, first for projection
estimators (\cite{LeuBar06}) and more recently for a more general
class of linear (affine) estimators (\cite{DalSal11}). Our further
analysis\vadjust{\goodbreak} will be based on this approach. Clearly, direct
sample splitting does not apply to independent samples that are
not identically distributed as in the present setup. Indeed, the
observations in the first sample no longer have the same
distribution as those in the second sample. On the other hand, the
approach based on Nemirovski's randomization can be still applied,
but it leads to somewhat weaker results involving an additional
expectation over a randomization distribution and a bigger remainder term than in our oracle inequalities.

\subsection{Aggregation of Linear Estimators}

Suppose that we are given a finite family $\{\hat f_1, \ldots, \allowbreak\hat
f_K\}$ of linear estimators defined by
\begin{equation}
\label{EQ:affine}
\hat f_j(x) = \bY^\top a_j(x),
\end{equation}
where $a_j(\cdot)$ are given functions with values in $\R^n$.  This
representation is quite general; for example, $\hat f_j$ can be
ordinary least squares, (kernel) ridge regression estimators or
diagonal linear filter estimators; see \citet{Kne94};
\citet{DalSal11} for a longer list of relevant examples. The vector
of values $(\hat f_j(x_i), i=1,\dots,n)$ equals to $A_j \bY$ where
$A_j$ an $n \times n$ matrix with rows $a_j(x_i), i=1,\dots,n$.

Now, we would like to consider mixtures of such estimators rather
than mixtures of deterministic functions as in the previous
sections. For this purpose, exponential weights have to be slightly
modified. Indeed, note that in Section~\ref{SEC:exp}, the risk of a
deterministic function $f_j$ is simply estimated by the empirical
risk $\hat R_n(f_j)$, which is plugged into the expression for the
weights. Clearly, $\E \hat R_n(f_j)= R(f_j)+\sigma^2$ so that $ \hat
R_n(f_j)$ is an \textit{unbiased} estimator of the risk $R(f_j)$ of
$f_j$ up to an additive constant. For a linear estimator $\hat f_j$
defined in~\eqref{EQ:affine}, $\hat R_n(\hat f_j)-\sigma^2$ is no
longer an unbiased estimator of the risk $\E R(\hat f_j)$. It is
well known that the risk of the linear estimator $\hat f_j$ has the
form\looseness=1
\[
\E R(\hat f_j)= \|(A_j -\bI)\eta\|^2 + \frac{\sigma^2}{n}
\Tr[A_j^\top A_j],
\]\looseness=0
where $\Tr[A]$ denotes the trace of a matrix $A$, and $\bI$ denotes
the $n \times n$ identity matrix. Moreover, an\vspace*{1pt} unbiased estimator of
$\E R(\hat f_j)$ is given by a version of Mallows's~$C_p$,
\begin{equation}
\label{EQ:unb}
\qquad\tilde R^{\mathrm{unb}}_n(\hat f_j)=\|\bY-\hat
f_j\|^2 + \frac{2\sigma^2}{n} \Tr[A_j] -\sigma^2.
\end{equation}
Then, for linear estimators, the exponential weights and the
corresponding aggregate\vadjust{\goodbreak} are modified as follows:
\begin{eqnarray}
\label{EQ:defexp2}
\hat \lambda^\pi_j &=&\frac{ \exp(-n\tilde
R^{\mathrm{unb}}_n(\hat f_j)/\beta)\pi_j}{\sum_{k=1}^K\exp(-n\tilde
R^{\mathrm{unb}}_n(\hat f_k)/\beta)\pi_k}, \nonumber\\[-8pt]\\[-8pt]
\hat f^\pi &=&
\sum_{k=1}^K \hat \lambda^\pi_k \hat f_k.\nonumber
\end{eqnarray}
Note that for deterministic $f_j$, we naturally define  $\tilde
R^{\mathrm{unb}}_n(f_j)=\hat R_n(f_j)-\sigma^2$, so that
definition~\eqref{EQ:defexp2} remains consistent
with~\eqref{EQ:defexp}. With this more general definition of
exponential weights, \citet{DalSal11} prove the following risk
bounds for the aggregate~$\hat f^\pi$.

\begin{TH1}
\label{TH:OIaff}
 Let $\{\hat f_1, \ldots, \hat f_K\}$ be a family of
linear estimators defined in~\eqref{EQ:affine} such that the
matrices $A_j$ are symmetric, positive definite and
$A_jA_k=A_kA_j$, 
for all $1\le j,k \le K$. Then the exponentially weighted\vspace*{1pt} aggregate
$\hat f^\pi$  defined in~\eqref{EQ:defexp2} with $\beta \ge
8\sigma^2$ satisfies
\begin{eqnarray}
\label{EQ:OIaff}
\hspace*{19pt}\E R(\hat f^\pi) &\le& \min_{\lambda \in \Lambda^K}
 \Biggl\{\sum_{j=1}^K \lambda_j \E R(\hat f_j) + \frac{\beta}{n}
\cK(\lambda, \pi)  \Biggr\},\hspace*{-19pt}
\\ \label{EQ:OIaff1}
\E R(\hat f^\pi) &\le&  \min_{j=1,\dots, K}  \biggl\{\E R(\hat f_j) +
\frac{\beta}{n} \log( \pi_j^{-1}) \biggr\}.
\end{eqnarray}
If all the $A_j$ are projection matrices
($A_j^\top =A_j$,\break \mbox{$A_j^2=A_j$}), then the above inequalities hold with\linebreak[4] $\beta \ge
4\sigma^2$.
\end{TH1}

Here, bound~\eqref{EQ:OIaff1} follows immediately
from~\eqref{EQ:OIaff}.  In the rest of the paper, we mainly use the
last part of this theorem concerning projection estimators. The
bound~\eqref{EQ:OIaff1} for this particular case was originally
proved in \citet{LeuBar06}.  The result of \citet{DalSal11} is, in
fact, more general than Theorem~\ref{TH:OIaff} covering nondiscrete
priors in the spirit of \eqref{EQ:PAC2}, and it applies not only to
linear, but also to affine estimators $\hat f_j$.

\section{Sparse Estimation}
\label{SEC:sparse}

The family of projection estimators that we consider in this section
is the family of all $2^M$ least squares estimators, each of which
is characterized by its sparsity pattern. We examine properties of
these estimators, and show that their mixtures with exponential
weights satisfy sparsity oracle inequalities for suitably chosen
priors $\pi$.

\subsection{Sparsity Pattern Aggregation}

Assume that we are given a dictionary of functions
$\cH=\{f_1, \ldots, f_M\}$. However,\vadjust{\goodbreak} we will not aggregate the
elements of the dictionary, but rather the least squares estimators
depending on all the $f_j$. We denote by $\bX$, the $n \times M$
\textit{design matrix} with elements $\bX_{i,j}=f_j(x_i)$, $i=1,
\ldots, n, j=1, \ldots, M$.

A \textit{sparsity pattern} is a binary vector $\pp \in
\cP:=\{0,1\}^M$. The terminology comes from the fact that the
coordinates $\pp_j$ of such vectors can be interpreted as indicators
of presence ($\pp_j=1$) or absence ($\pp_j=0$) of a given feature
indexed by $j\in\{1,\ldots,M\}$. We denote by $|\pp|$ the number of
ones in the sparsity pattern $\pp$, and by $S^{\pp}$ the linear span
of canonical basis vectors $e^{(j)}$, such that $\pp_j=1$.

For $\pp \in \cP$, let $\hat \theta_\pp$ be any least squares
estimator on~$S^{\pp}$ defined by
\begin{equation}
\label{EQ:lse}
\qquad\hat \theta_{\pp} \in \argmin_{\theta \in S^{\pp}}
\|\bY-\ff_\theta\|^2  
\quad  \mbox{with }   \ff_\theta=\sum_{j=1}^M \theta_jf_j.
\end{equation}
The following simple lemma gives an oracle inequality for the least
squares estimator.  It follows easily from the  Pythagorean theorem.
Moreover, the random variables $\xi_1, \ldots, \xi_n$ need not be
Gaussian for the result to hold.

\begin{lem}\label{LEM:lse}
Fix $\pp \in \cP$.  Then any least squares estimator $\hat \theta_\pp$
defined in \eqref{EQ:lse} satisfies
\begin{eqnarray}
\label{EQ:pytha}
\E\|\ff_{\hat \theta_\pp}-\eta\|^2&=& \min_{\theta
\in S^{\pp} }\|\ff_{ \theta}-\eta\|^2 + \sigma^2\frac{d_{\pp}}{n}\nonumber\\[-8pt]\\[-8pt]
&\le& \min_{\theta \in S^{\pp} }\|\ff_{ \theta}-\eta\|^2 +
\sigma^2\frac{|\pp|}{n},\nonumber
\end{eqnarray}
where $d_\pp$ is the dimension of the linear subspace
$\{\bX\theta \dvtx  \theta \in S^{\pp} \}$ .
\end{lem}

Clearly, if $|\pp|$ is small compared to $n$, the oracle inequality
gives a good performance guarantee for the least squares aggregate
$\ff_{\hat \theta_\pp}$. Nevertheless, it may be the case that the
approximation error $\min_{\theta \in S^{\pp}}\|\ff_{
\theta}-\eta\|^2 $ is quite large. Hence, we are looking for a
sparsity pattern such that $|\pp|$ is small and that yields a least
squares aggregate with small approximation error. This is clearly a
model selection problem, as described in Section~\ref{sec1}.

Observe that for each sparsity pattern $\pp \in \cP$, the function
$\ff_{\hat \theta_\pp}$ is a projection estimator of the form
$\ff_{\hat \theta_\pp}=A_\pp \bY$ where the $n\times n$ matrix
$A_\pp$ is the projector onto $\{\bX\theta \dvtx  \theta \in S^{\pp}
\}$ 
(as above, we
identify the functions $f_j, \ff_{\hat \theta_\pp}$ with the vectors
of their values at points $x_1,\dots,x_n$ since the risk depends
only on these values). Therefore $\Tr[A_\pp]=d_{\pp}$. We have seen
in the previous section that, to solve the problem of model
selection, projection\vadjust{\goodbreak} estimators can be aggregated using exponential
weights. Thus, instead of selecting the best sparsity pattern, we
resort to taking convex combinations leading to what is called
\textit{sparsity pattern aggregation}. For any sparsity pattern $\pp
\in \cP$, define the exponential weights $\hat \lambda^\pi_\pp$ and
the \textit{sparsity pattern aggregate}~$\tilde f^{\pi}$, respectively,
by
\begin{eqnarray}
\hat \lambda^\pi_\pp &=& \frac{\exp(-n\tilde
R^{\mathrm{unb}}_n(\ff_{\hat
\theta_\pp})/\beta)\pi_\pp}{\sum_{\pp'\in\cP}\exp(-n\tilde
R^{\mathrm{unb}}_n(\ff_{\hat \theta_{\pp'}})/\beta)\pi_{\pp'}},\nonumber\\
 \tilde f^{\pi} &=& \sum_{\pp \in \cP} \hat \lambda^\pi_\pp
\ff_{\hat \theta_\pp},\nonumber
\end{eqnarray}
where $\pi=(\pi_{\pp})_{\pp \in \cP}$ is a probability distribution
(prior) on the set of sparsity patterns $\cP$.


To study the performance of this method, we can now apply the last
part of Theorem~\ref{TH:OIaff} dealing with projection matrices. Let
$\pp(\theta) \in \cP$ be the sparsity pattern of $\theta\in \R^M$,
that is, a vector with components $\pp_j(\theta)=1$ if $\theta_j \neq
0$, and $\pp_j(\theta)=0$ otherwise. Note that
$|\pp(\theta)|=|\theta|_0$. Combining~\eqref{EQ:OIaff1} and
Lemma~\ref{LEM:lse} and the fact that $\{\theta\dvtx
\pp(\theta)=\pp\}\subset S^{\pp}$, we get that for $\beta\ge 4\sigma^2$
\begin{eqnarray}\label{EQ:OIspa}
\E R(\tilde f^{\pi}) &\le& \min_{\pp \in \cP}  \biggl\{ \E R(\ff_{\hat
\theta_\pp}) + \frac{\beta}{n}
\log( \pi_\pp^{-1})  \biggr\} \nonumber
\\
&\le & \min_{\pp \in \cP} \biggl\{ \min_{\theta\dvtx
\pp(\theta)=\pp}\|\ff_{ \theta}-\eta\|^2+ \sigma^2\frac{|\pp|}{n}\nonumber\\
&&\hspace*{77pt}{}  +
\frac{\beta}{n} \log( \pi_{\pp}^{-1})  \biggr\}
\\
&= & \min_{\theta \in \R^M} \biggl\{ \|\ff_{ \theta}-\eta\|^2 +
\sigma^2\frac{|\theta|_0}{n}\nonumber\\
&&\hspace*{48pt}{} + \frac{\beta}{n} \log\bigl(
\pi_{\pp(\theta)}^{-1}\bigr)  \biggr\},   \nonumber
\end{eqnarray}
where we have used that $\min_{\theta \in \R^M}$ can be represented
as $\min_{\pp \in \cP}\min_{\theta\dvtx   \pp(\theta)=\pp}$.

The remainder term in the balanced oracle
inequality~\eqref{EQ:OIspa} depends on the choice of the prior
$\pi$. Several choices can be considered depending on the
information that we have about the oracle, that is, about a potentially
good candidate $\theta$ that we would like to mimic. For example, we
can assume that there exists a good $\theta$ that is coordinatewise
sparse, group sparse or even that $\theta$ is piecewise constant.
While this approach to structure the prior knowledge seems to fit in
a Bayesian framework, we only pursue a frequentist setup. Indeed,
our risk measure is not averaged over a prior. Such priors on good
candidates for estimation are often used in a non-Bayesian
framework. For\vadjust{\goodbreak} example, in nonparametric estimation, it is usually
assumed that a good candidate function is smooth. Without such
assumptions, one may face difficulties in performing meaningful
theoretical analysis.

\subsection{Sparsity Priors}
\label{sub:priors}

\subsubsection{Coordinatewise sparsity}

This is the basic\break and most commonly used form of sparsity. The prior
$\pi$ should favor vectors $\theta$ that have a small number\break of
nonzero coordinates. Several priors have been\break \mbox{suggested} for this
purpose;
cf. \citet{LeuBar06}; \citet{Gir07}; \citet{RigTsy11}; \citet{AlqLou11}.
We consider here yet another prior, close to that of \citet{Gir07}.
The main difference is that the prior $\pi^{\textsc{c}}$ below
exponentially downweights sparsity patterns with large $|\pp|$,
whereas the prior in \citet{Gir08} downweights such patterns
polynomially. Define
\[
\pi^{\textsc{c}}_\pp= \Biggl[{M \choose
|\pp|}e^{|\pp|}H_M \Biggr]^{-1}, \quad H_M=\sum_{k=0}^M e^{-k}\le
\frac{e}{e-1}.
\]
It can be easily seen that $\sum_{\pp \in
\cP}\pi^{\textsc{c}}_\pp=1$ so that
$\pi^{\textsc{c}}=(\pi_\pp^{\textsc{c}},\pp\in\cP)$ is a probability
measure on $\cP$. Note that
\begin{eqnarray}
\label{EQ:bound_pi_coord}
\log [  ( \pi^{\textsc{c}}_\pp
 )^{-1} ]&=&\log {M \choose |\pp|} + |\pp| +\log (H_M)\nonumber\\[-8pt]\\[-8pt]
 & \le& 2|\pp|\log
 \biggl(\frac{eM}{|\pp|} \biggr) + \frac{1}{2},\nonumber
\end{eqnarray}
where we have used the inequality ${M \choose |\pp|} \le
 (\frac{eM}{|\pp|} )^{|\pp|}$ for $|\pp|\ne0$ and the
convention $0\log(\infty)=0$ for $|\pp|=0$.
Define the sparsity pattern aggregate
\begin{equation}
\label{EQ:spaC}
\tilde f^{\textsc{c}}= \sum_{\pp \in \cP} \hat \lambda^{\pi^{\textsc{c}}}_\pp\ff_{\hat \theta_\pp},
\end{equation}
where $\lambda^{\pi^{\textsc{c}}}_\pp$ is the exponential weight given
in~\eqref{EQ:defexp2}, and $\hat  \theta_\pp$ is the least squares
estimator~\eqref{EQ:lse}.

Plugging~\eqref{EQ:bound_pi_coord}  into~\eqref{EQ:OIspa} with
$\pi=\pi^{\textsc{c}}$ and $\beta = 4\sigma^2$ yields the following
sparsity oracle inequality:
\begin{eqnarray}
\label{EQ:spaC_final}
\qquad\E R(\tilde f^{{\textsc{c}}})&\le& \inf_{\theta
\in \R^M} \biggl\{ \|\ff_{ \theta}-\eta\|^2   \nonumber\\[-8pt]\\[-8pt]
 &&{}+
 \frac{9\sigma^2}{n}|\theta|_0 \log  \biggl(\frac{eM}{|\theta|_0} \biggr) +\frac{2 \sigma^2}{n}  \biggr\}.\nonumber
\end{eqnarray}
It is important to note that \eqref{EQ:spaC_final} is valid under
\textit{no assumption in the dictionary}. This is in contrast to the
Lasso and assimilated penalized procedures that are known to have
similar properties only under strong conditions on $\mathbf{X}$, such
as restricted isometry or restricted eigenvalue
conditions\vadjust{\goodbreak} (see, e.g., \cite{CanTao07}; \cite{BicRitTsy09};\linebreak[4]  \cite{KolLouTsy11}).

Another choice for $\pi$ in the framework of coordinatewise sparsity
can be found in \citet{RigTsy11} and yields the \textit{exponential
screening} estimator. The exponential screening aggregate satisfies
an improved version of the above sparsity oracle inequality with $|\theta|_0$ replaced by
$\min(|\theta|_0, \D)$ where $\D$ is the rank of the design
matrix~$\bX$. In particular, if the rank $\D$ is small, the
exponential screening aggregate adapts to it. Moreover, it is shown
in \citet{RigTsy11} that the remainder term of the oracle inequality
is optimal in a minimax sense.
%

\subsubsection{Fused sparsity}

When there exists a natural order among the functions $f_1, \ldots,
f_M$ in the dictionary, it may be appropriate to assume that there
exists a ``piecewise constant'' $\theta \in \R^M$, that is, $\theta$
with components taking only a small number of values, such that
$\ff_\theta$ has good approximation properties. This property often
referred to as \textit{fused sparsity} has been exploited in the image
denoising literature for two decades, originating with the classical
paper by \citet{RudOshFat92}. The \textit{fused Lasso} was
introduced in \citet{TibSauRos05} to deal with the same problem in
one dimension instead of two. Here we suggest another method that
takes advantage of fused sparsity using the idea of mixing with
exponential weights. Its theoretical advantages are demonstrated by
the sparsity oracle inequality in Corollary~\ref{COR:fus} below.

At first sight, this problem appears to be different from the one
considered above since a good $\theta \in \R^M$ need not be sparse.
Yet, the fused sparsity assumption on $\theta$ can be reformulated
into a coordinatewise sparsity assumption. Indeed, let $D$ be the
$M\times M$ matrix defined by the relations $(D\theta)_1=\theta_1$
and $(D\theta)_j=\theta_j-\theta_{j-1}$ for $j=2, \ldots, M$, where
$(D\theta)_j$ is the $j$th component of~$D\theta$. We will call $D$
the ``first differences'' matrix. Then $\theta$ is fused sparse if
$|D\theta|_0$ is small.

We now consider a more general setting with an arbitrary invertible
matrix~$D$, again declaring $\theta$ to be fused sparse if
$|D\theta|_0$ is small. Possible definitions of $D$ can be based on
higher order differences or combinations of differences of several
orders accounting for other types of sparsity.
 For each sparsity pattern
$\pp \in \cP$, we define the least squares estimator
\begin{equation}
\label{EQ:lse-D}
\hat \theta_{\pp}^D \in \argmin_{\theta \in
\R^M\dvtx D\theta\in S^{\pp}} \|\bY-\ff_\theta\|^2.
\end{equation}
The corresponding estimator $\ff_{\hat \theta_{\pp}^D}$ (as
previously, without loss of generality we consider $\ff_{\hat
\theta_{\pp}^D}$ as an $n$-vector) takes the form $\ff_{\hat
\theta_{\pp}^D}=A_\pp^D\bY$ where $A_\pp^D$ is the projector onto
the linear space $\cL_\pp=\{\bX\theta \dvtx  \theta \in
\R^M, D\theta\in S^{\pp}\}$. In particular,
$\Tr[A_\pp^D]=\mathrm{dim}(\cL_{\pp})$, where
$\mathrm{dim}(\cL_{\pp})$ is the dimension of $\cL_\pp$. Moreover,
it is straightforward to obtain the following result, analogous to
Lem\-ma~\ref{LEM:lse}.

\begin{lem}\label{LEM:lse-D}
Fix $\pp \in \cP$, and let $D$ be an invertible matrix.  Then any least squares estimator $\hat \theta_\pp^D$
defined in \eqref{EQ:lse-D} 
satisfies
\begin{eqnarray}
\label{EQ:pytha-D}
\qquad\E\|\ff_{\hat \theta_\pp^D}-\eta\|^2&=&
\mathop{\min_{\theta \in \R^M:}}_{D\theta\in S^{\pp}}
\|\ff_{
\theta}-\eta\|^2 + \sigma^2\frac{\mathrm{dim}(\cL_{\pp})}{n}\nonumber\\[-8pt]\\[-8pt]
 &\le&
\mathop{\min_{\theta \in \R^M:}}_{D\theta\in S^{\pp}}
\|\ff_{
\theta}-\eta\|^2 + \sigma^2\frac{|\pp|}{n}.\nonumber
\end{eqnarray}
\end{lem}

We are therefore in a position to apply the results from
Section~\ref{SEC:estim}. For example, if $\pp$ is sparse, and $D$ is
the ``first differences'' matrix, the least squares estimator $\hat
\theta_{\pp}^D$ is piecewise constant with a small number $|\pp|$ of
jumps.

Now, since the problem has been reduced to coordinatewise sparsity,
we can choose the prior $\pi^{\textsc{c}}$ to favor vectors $\theta
\in \R^M$ that are piecewise constant with a small number of jumps.
Define the fused sparsity\vspace*{1pt} pattern aggregate $\tilde f^{\textsc{f}}$
by
\begin{equation}
\label{EQ:spaF}
\tilde f^{\textsc{f}}=\sum_{\pp \in \cP} \hat \lambda^{\pi^{\textsc{c}}}_\pp\ff_{\hat \theta_\pp^D},
\end{equation}
where $\hat \lambda^{\pi^{\textsc{c}}}_\pp$ is the exponential weight
defined in~\eqref{EQ:defexp2}, and $\hat  \theta_\pp^D$ is the least
squares estimator defined in~\eqref{EQ:lse-D}. Note that we can combine
\eqref{EQ:OIaff1} with Lemma~\ref{LEM:lse-D} in the same way as in
\eqref{EQ:OIspa} with the only difference that we use now the relation
$\min_{\pp \in \cP}\min_{\theta \dvtx \pp(D\theta)=\pp}(\cdot)=
\min_{\theta \dvtx   D\theta \in \R^M}(\cdot)=\min_{\theta \in
\R^M}(\cdot)$. This and~\eqref{EQ:bound_pi_coord} imply the following
bound.

\begin{cor}\label{COR:fus}
Let $D$ be an invertible matrix. The fused sparsity pattern aggregate $\tilde f^{\textsc{f}}$
defined in \eqref{EQ:spaF} with $\beta = 4\sigma^2$ satisfies
\begin{eqnarray}
\label{EQ:spaF1}
&&\hspace*{20pt}\E R(\tilde f^{{\textsc{f}}}) \nonumber\\
 &&\hspace*{20pt}\quad\le \inf_{\theta \in
\R^M} \biggl\{ \|\ff_{ \theta}-\eta\|^2\\
&&\hspace*{20pt}\quad\hspace*{38pt}{}  +
 \frac{9\sigma^2}{n}|D\theta|_0 \log
  \biggl(\frac{eM}{|D\theta|_0} \biggr)+\frac{2 \sigma^2}{n}   \biggr\}.\nonumber
\end{eqnarray}
\end{cor}

 To our knowledge, analogous bounds for fused\break
Lasso are not available.\vadjust{\goodbreak} Furthermore, Corollary~\ref{COR:fus} holds
under no assumption on the matrix $\mathbf{X}$, which cannot be the
case for the Lasso type methods. Let us also emphasize that
Corollary~\ref{COR:fus} is valid for any invertible matrix $D$,  and
not only for the standard ``first differences'' matrix $D$ defined above.

\subsubsection{Group sparsity}

Since recently, estimation under group sparsity has been intensively
discussed in the literature. Starting from \citet{YuaLin06},
several estimators have been studied, essentially the Group Lasso
and some related penalized techniques. Theoretical properties of the
Group Lasso are treated in some generality by \citet{HuaZha10} and \citet{LouPonTsy11}
where one can find further references.
Here we show that one can deal with group sparsity using
exponentially weighted aggregates. The new estimator that we propose
presents some theoretical advantages as compared to the\break Group Lasso
type methods.

Let $B_1,\ldots, B_K$ be given subsets of $\{1,\ldots, M\}$ cal\-led the
groups. We impose no restriction on $B_j$'s; for example, they need not
form a partition of $\{1,\ldots, M\}$ and can overlap. In this section,
we consider $\theta\in\R^M$ such that
$\operatorname{supp}(\theta)\subseteq B \triangleq\bigcup_{k=1}^K B_k $
where $\operatorname{supp}(\theta)$ is the support of $\theta$. For any
such $\theta$, we denote by $J(\theta)$ the subset of $\{1,\ldots, K\}$
of smallest cardinality among all $J$ satisfying
$\operatorname{supp}(\theta)\subseteq B_J \triangleq\bigcup_{k\in J}
B_k$. We assume without loss of generality that $J(\theta)$ is unique.
(If there are several $J$ of same cardinality satisfying this property,
we define $J(\theta)$ as the smallest among them with respect to some
partial ordering of the subsets of $\{1,\ldots, K\}$.) Set
\[
g(\theta)=|J(\theta)|,\quad B(\theta)=\bigcup_{k\in J(\theta)}B_k.
\]
The group sparsity setup assumes that there exists $\theta\in\R^M$
such that $\|\ff_\theta-\eta\|^2$ is small and that $\theta$ is
supported by a small number of groups, that is, that
$g(\theta)\ll K$.

Let now $J$ be a subset $\{1,\ldots, K\}$. Denote by $\pp^{J}$ the
sparsity pattern with coordinates defined by
\[
\pp^{J}_j=\cases{
1, & if $j\in B_J$,\cr
0, & otherwise,}
\]
for $j=1,\ldots,M$. Consider the set of all such sparsity patterns:
\[
\cP_{\mathrm{G}}=\bigl\{\pp^{J}, J \subseteq\{1,\ldots,
K\}\bigr\}.
\]
To each sparsity pattern $\pp^J \in\cP_{\mathrm{G}}$ we assign a least
squares estimator $\hat\theta_{\pp^J}$, cf. \eqref{EQ:lse},
constrained to having null coordinates outside of $B_J=\bigcup_{k \in
J} B_k$.

Define the following prior on $\cP_{\mathrm{G}}$:
\[
\pi^{\mathrm{G}}_{\pp^J}= \left[\pmatrix{K \cr
|J|}e^{|J|}H_K\right]^{-1},\quad J \subseteq\{1,\ldots, K\}.
\]
This prior enforces group sparsity by favoring the small number of
groups $|J|$. As in \eqref{EQ:bound_pi_coord}, we obtain
\begin{equation}
\label{EQ:bound_pi_group} \log[ ( \pi^{\mathrm{G}}_{\pp^J}
)^{-1}] \le2|J|\log\biggl(\frac{eK}{|J|}\biggr) + \frac
{1}{2}.
\end{equation}
We introduce now the sparsity pattern aggregate
\begin{equation}
\label{EQ:spaG}
\tilde f^{\mathrm{G}}= \sum_{\pp\in\cP_{\mathrm{G}}}
\hat\lambda^{\pi^{\mathrm{G}}}_\pp\ff_{\hat\theta_\pp},
\end{equation}
where $\lambda^{\pi^{\mathrm{G}}}_\pp$ is the exponential weight
defined in \eqref{EQ:defexp2} and $\hat \theta_\pp$ is the least
squares estimator defined in \eqref{EQ:lse}.

For any $\pp=\pp^J\in\cP_{\mathrm{G}}$, we have
$V_J\triangleq\{\theta\dvtx \operatorname{supp}(\theta)\subseteq B,
J(\theta)=J\}\subseteq S^{\pp}$. Arguing as in \eqref{EQ:OIspa} with
$\cP_{\mathrm{G}}$ instead of $\cP$, setting $\pi=\pi^{\mathrm{G}}$ and
using \eqref{EQ:bound_pi_group} we obtain that, for $\beta\ge
4\sigma^2$,
\begin{eqnarray*}
&&
\E R (\tilde f^{{\mathrm{G}}})\nonumber\\
&&\quad\le \min_{\pp\in\cP_{\mathrm{G}}}
\biggl\{ \min_{\theta\in S^{\pp}}\|\ff_{ \theta}-\eta\|^2 +
\sigma^2\frac{|\pp|}{n} \nonumber\\
&&\hspace*{48.3pt}\qquad{}+
\frac{\beta}{n} \log( (\pi_{\pp}^{\mathrm{G}})^{-1}) \biggr\}\nonumber
\nonumber\\
&&\quad\le \min_{J \subseteq\{1,\ldots,K\}}
\biggl\{ \min_{ \theta\in V_J}\|\ff_{ \theta}-\eta\|^2 +
\sigma^2\frac{|B_J|}{n} \\
&&\qquad\hspace*{76.4pt}{}+
\frac{\beta}{n} \log( (\pi_{\pp^J}^{\mathrm{G}})^{-1}) \biggr\}
\nonumber\\
&&\quad\le \min_{J \subseteq\{1,\ldots,K\}}
\min_{ \theta\in V_J}
\biggl\{\|\ff_{ \theta}-\eta\|^2 +
\sigma^2\frac{|B(\theta)|}{n} \nonumber\\
&&\qquad\hspace*{67.5pt}{}+ \frac{2\beta}{n} g(\theta)\log
\biggl(\frac{eK}{g(\theta)}\biggr) + \frac{\beta}{2n}\biggr\}.
\nonumber
\end{eqnarray*}
This leads to the following oracle inequality.

\begin{cor}\label{COR:group}
The group sparsity pattern aggregate $\tilde f^{\mathrm{G}}$ defined in
\eqref{EQ:spaG} with $\beta= 4\sigma^2$ satisfies
\begin{eqnarray}\label{EQ:COR:group}
\E R(\tilde f^{{\mathrm{G}}})&\le&
\mathop{\inf_{\theta\in\R^M:}}_{
\operatorname{supp}(\theta)\subseteq B }
\biggl\{ \|\ff_{ \theta}-\eta\|^2
+ \sigma^2\frac{|B(\theta)|+2}{n}\nonumber\\[-8pt]\\[-8pt]
&&\qquad\hspace*{35pt}{} + \frac{8\sigma^2}{n}g(\theta)
\log\biggl(\frac{eK}{g(\theta)}\biggr) \biggr\}.\hspace*{-10pt}\nonumber
\end{eqnarray}
\end{cor}

We see from Corollary~\ref{COR:group} that if there exists an ideal
``oracle'' $\theta$ in $\R^M$, such that the approximation error $
\|\ff_{ \theta}-\eta\|^2$ is small, and $\theta$ is sparse in the
sense that it is supported by a small number of groups, then the
sparsity pattern aggregate $\tilde f^{{\textsc{g}}}$ mimics the risk
of this oracle.\vadjust{\goodbreak}

A remarkable fact is that Corollary~\ref{COR:group} holds for
arbitrary choice of groups~$B_j$. They can overlap and not
necessarily cover the whole set $\{1,\dots,M\}$.

To illustrate the power of the oracle
inequality~\eqref{EQ:COR:group}, we consider the multi-task learning
setup as in \citet{LouPonTsy11}. Namely, assume that all the groups
$B_j$ are of the same size $T$ and form a partition of
$\{1,\dots,M\}$, so that $M=KT$. We restrict our analysis to the
class $\cF_s$ of regression functions $\eta$ such that
$\eta=\ff_{\theta}$ for some $\theta$ satisfying $g(\theta)\le s$
where $s\le K$ is a given integer. Then $|B(\theta)|\le sT$.
Combining these remarks with~\eqref{EQ:COR:group} and with the fact
that the function $x\mapsto x\log  (\frac{eK}{x} )$ is
increasing, we find that, uniformly over $\eta\in\cF_s$,
\begin{equation}\label{EQ:COR:group1}
\qquad\E R(\tilde f^{{\textsc{g}}})\le   \frac{\sigma^2s}{n} \biggl(T +
8\log  \biggl(\frac{eK}{s} \biggr)+\frac{2}{s} \biggr).
\end{equation}
On the other hand, a minimax lower bound on the same class $\cF_s$
is available in \citet{LouPonTsy11}. It has exactly the form of the
right-hand side of~\eqref{EQ:COR:group1}; cf. equation (6.2)
in \citet{LouPonTsy11}. This immediately implies that (i) the lower
bound of \citet{LouPonTsy11} is tight so that $\frac{s}{n} (T +
\log  (\frac{K}{s} ) )$
is the optimal rate of
convergence on $\cF_s$, and (ii) the estimator $\tilde
f^{{\textsc{g}}}$ is rate optimal. To our knowledge, this gives the
first example of rate optimal estimator under group sparsity. The
upper bounds for the Group Lasso estimators in
\citet{HuaZha10} and \citet{LouPonTsy11} as well as in
the earlier papers cited therein depart from this optimal rate at
least by a
 logarithmic factor. Furthermore, they are obtained under
strong assumptions on the dictionary such as restricted isometry or
restricted eigenvalue type conditions, while~\eqref{EQ:COR:group1}
is valid under no assumption on the dictionary.

\section{Related Problems}\label{SEC:related}

In this paper, we have considered only the Gaussian regression model
with fixed design and known variance of the noise. This is a basic
setup where the sharpest results, expressed in terms of
sparsity oracle inequalities, are now available for exponentially
weighted (EW) procedures both in aggregation and sparsity scenarios.
Similar but somewhat weaker properties are obtained for exponential weighting
in several other models.

\textit{Models with i.i.d. observations}. Some EW aggregates achieve
sparsity oracle inequalities in regression model with random
design (\cite{DalTsy10}; \cite{AlqLou11}; \cite{Ger11}) as well as in density
estimation and classification problems (\citeauthor{DalTsy10}\break (\citeyear{DalTsy10})). However,
the results differ in several aspects from those of the present
paper. First, they do not use aggregation of estimators, but rather
EW procedures driven by continuous priors (\cite{DalTsy10}; \cite{Ger11}), or
by priors with both discrete and continuous components
(\cite{AlqLou11}). The developments in \citet{DalTsy10}; \citet{AlqLou11}; \citet{Ger11} start from the general oracle inequalities similar to
(\ref{EQ:PAC2}), which are sometimes called PAC-bounds; cf. recent
overview in\break \citet{Cat07}. Sparsity oracle inequalities are then
derived from PAC-bounds. However, as opposed\break to~(\ref{EQ:spaC_final}), they involve not only $|\theta|_0$ but also
the $\ell_1$-norm of~$\theta$. The estimators in \citet{DalTsy10}; \citet{Ger11} are defined as an average of exponentially weighted
aggregates over the sample sizes from 1 to $n$. This is related to
earlier work on mirror averaging; cf. \citet{JudNazTsy05}; \citet{JudRigTsy08}, which in turn, is inspired by the concept of mirror
descent in optimization due to Nemirovski. Finally, the
computational algorithms are also quite different from those that we
describe in the next section. For example, under continuous sparsity
priors, one of the suggestions is to use Langevin Monte-Carlo; cf.
Dalalyan and Tsybakov (\citeyear{DalTsy12,DalTsy10}).

\textit{Unknown variance of the noise, non-Gaussian\break noise}.
Modifications of EW procedures and of the corresponding oracle
inequalities for the case of unknown variance $\sigma^2$ are
discussed in \citet{Gir07}; \citet{Ger11}. Moreover, the results can be
extended to regression with non-Gaussian noise under deterministic
or random design (\cite{DalTsy07}; \cite{DalTsy08}; \cite{DalTsy12}; \cite{Ger11}). In
particular, \cite{Ger11} uses a version of the EW estimator with
data-driven truncation to cover a rather general  noise structure. The
estimator satisfies a balanced oracle inequality but not a sparsity
oracle inequality as defined here, since along with $|\theta|_0$, it
involves other characteristics of $\theta$ and of the target
function $\eta$.

\section{Numerical Implementation}
\label{SEC:num}

\label{SEC:simul}


All the sparsity pattern aggregates defined in the previous section
are of the form $\ff_{\theta^{\exp}}$, where
\begin{equation}
\label{EQ:theta-exp}
\theta^{\exp} =\sum_{\pp \in \cG} \hat
\lambda_\pp^\pi \bar \theta_\pp
\end{equation}
for some $\cG \subset \cP$, $\lambda_\pp^\pi$ is the exponential
weight defined in~\eqref{EQ:defexp2}, and $\bar \theta_\pp$ is
either $\hat \theta_\pp$ defined in~\eqref{EQ:lse} or $\hat
\theta_\pp^D$ defined in~\eqref{EQ:lse-D}.

From~\eqref{EQ:theta-exp}, it is clear that one needs to add up with
some weights $2^M$ (or $2^K$ in the case of group sparsity with $K$
groups) least squares estimators to compute $\theta^{\exp}$
\textit{exactly}. In many applications this number is prohibitively
large. However, most of the terms in the sum receive an
exponentially low weight with the choices of $\pi$ that we have
described. We resort to a numerical approximation that exploits this
fact.

\begin{figure}
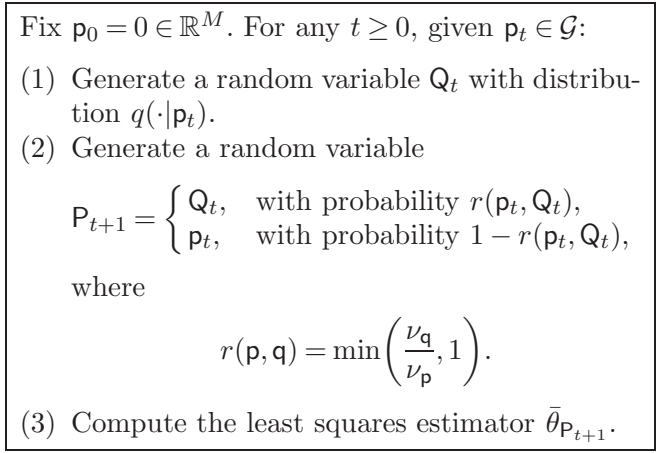

\bookbox{
Fix $\pp_0=0\in \R^M$. For any $t\ge 0$, given $\pp_t \in \cG$:
\begin{enumerate}[(3)]
\item[(1)] Generate a random variable $\qQ_{t}$ with distribution
$q(\cdot|\pp_t)$.
\item[(2)] Generate a random variable
\[
\pP_{t+1}= \cases{
\qQ_t, & with probability $r(\pp_t, \qQ_t)$,\cr
\pp_t, & with probability $1-r(\pp_t, \qQ_t)$,
}
\]
where
\[
r(\pp, \qqq)=\min \biggl( \frac{\nu_{\qqq}}{\nu_{\pp}},1 \biggr).
\]
\item[(3)] Compute the least squares estimator $\bar \theta_{\pP_{t+1}}$.
\end{enumerate}}
\caption{The Metropolis--Hastings algorithm on the $M$-hypercube.}\label{FIG:MH}
\end{figure}

Note that $\theta^{\exp}$ is obtained
as the expectation of the random variable $\hat \theta_{\sf P}$
or $\hat \theta_\pp^D$ where ${\sf P}$ is a random
variable taking values in $\cP$ with probability distribution $\nu$
given by
\[
\nu_\pp= \frac{\exp(-n\tilde R^{\mathrm{unb}}_n(\ff_{\bar
\theta_\pp})/\beta)\pi_\pp}{\sum_{\pp'\in \cG}\exp(-n\tilde
R^{\mathrm{unb}}_n(\ff_{\bar \theta_{\pp'}})/\beta)\pi_{\pp'}},
\quad  \pp \in \cG \subset \cP.
\]
This Gibbs-type distribution can be expressed as the stationary
distribution of the Markov chain generated by the
Metropolis--Hastings (MH) algorithm (see, e.g., \cite{RobCas04}, Section~7.3). We now describe the MH algorithm
employed here. Note that in the examples considered in the previous
section, $\mathcal{G}$ is either the hypercube $\cP$ or the
hypercube $\cP_{\textsc{g}}$. For any $\pp \in \mathcal{G}$, define
the instrumental distribution $q(\cdot| \pp)$ as the uniform
distribution on the neighbors of $\pp$ in $\mathcal{G}$, and notice
that since each vertex has the same number of neighbors, we have
$q(\pp| \qqq)=q(\qqq| \pp)$ for any $\pp, \qqq \in \cP$. The MH
algorithm is defined in Figure~\ref{FIG:MH}. We use here the uniform
instrumental distribution for the sake of simplicity. Our
simulations show that it yields satisfactory results both in
terms of performance and of speed. Another choice of $q(\cdot|\cdot)$ can
potentially further accelerate the convergence of the MH algorithm.

From the results of \citet{RobCas04} (see also \cite{RigTsy11},
Theorem~7.1) the Markov chain $(\pP_t)_{t\ge 0}$ defined in
Figure~\ref{FIG:MH} is ergodic. In other words, it holds
\[
\lim_{T\to \infty}\frac{1}{T}\sum_{t=T_0+1}^{T_0+T}\bar \theta_{\pP_t}
=\sum_{\pp \in \cG}\bar \theta_\pp \nu_\pp,\quad   \mbox{almost surely},
\]
where $T_0\ge 0$ is an arbitrary integer.

In view of this result, we approximate $\theta^{\exp}=\break\sum_{\pp \in \cG}\bar \theta_\pp \nu_\pp$ by
\[
\tilde{\theta}^{\exp}_T=
\frac{1}{T}\sum_{t=T_0+1}^{T_0+T}\bar
\theta_{\pP_t},
\]
which is close to ${\theta^{\exp}}$ for sufficiently large $T$. One
remarkable feature of the MH algorithm is that it involves only the
ratios $\nu_\qqq/\nu_\pp$ where $\pp$ and $\qqq$ are two neighbors
in~$\mathcal{G}$. Such ratios are easy to compute, at least in the
examples given in the previous section. As a result, the MH
algorithm in this case takes the form of a stochastic greedy
algorithm with averaging, which measures a trade-off between sparsity
and prediction to decide whether to add or remove a variable. In all
subsequent examples, we use a pure R implementation of the sparsity
pattern aggregates. While the benchmark estimators considered below
employ a C based code optimized for speed, we observed that a safe
implementation of the MH algorithm (three times more iterations than
needed) exhibits an increase of computation time of at most a factor
two.

\subsection{Numerical Experiments}

The aim of this subsection is to illustrate the performance of the
sparsity pattern aggregates $\tilde f^{\textsc{c}}$ and $\tilde
f^{\textsc{f}}$ defined in~\eqref{EQ:spaC} and~\eqref{EQ:spaF}
respectively, on a simulated dataset and to compare it with
state-of-the-art procedures in sparse estimation. In our
implementation, we replace the prior $\pi^{\textsc{c}}$ by the
\textit{exponential screening} prior employed in \citet{RigTsy11}. As
a result, the following results are about the exponential screening
(\textsc{es}) aggregate defined in \citet{RigTsy11}. Nevertheless, it
presents the same qualitative behavior as the aggregates constructed
above.

\begin{figure*}[t]

\includegraphics{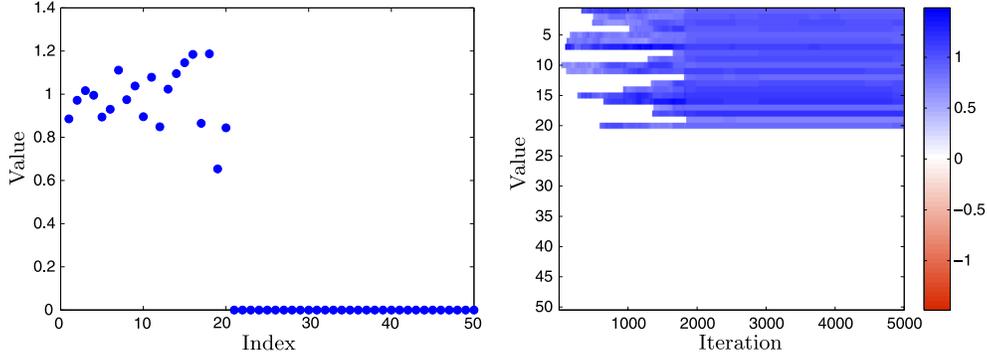}

\caption{Typical realization for $(M,n,S)=(500,200,20)$. \textit{Left:} Value of the  $\tilde{\theta}^{\exp}_T$, $T_0=3000,  T=7000$.
\textit{Right:} Value of $\hat \theta_{\pP_t}$ for $t=1, \ldots, 5000$. Only the first 50 coordinates are shown for each vector.}\vspace*{-5pt}
\label{FIG:iter}\vspace*{-3pt}
\end{figure*}

While our results for the \textsc{es} estimator hold under no
assumption on the dictionary, we compare the behavior of our
algorithm in a well-known example where sparse estimation by
$\ell_1$-penalized techniques is theoretically achievable.

Consider the model $\bY=\bX\theta^*+\sigma \xi$, where $\bX$ is an
$n\times M$ matrix with independent standard Gaussian
entries, and $\xi \in \R^n$ is a vector of independent standard Gaussian
random variables and is independent of $\bX$. Depending on our sparsity assumption, we choose two different $\theta^*$.

The variance is chosen as
$\sigma^2=\|\ff_{\theta^*}\|^2/9=|\bX{\theta^*}|_2^2/\allowbreak(9n)$ following
the numerical experiments of Candes and Tao [(\citeyear{CanTao07}), Section~4]. Here
$|\cdot|_2$ denotes the $\ell_2$ norm. For different values of
$(n,M,S)$, we run the \textsc{es} algorithm on 500 replications of the
problem and compare our results with several other popular
estimators in the literature on sparse estimation that are readily
implemented in R. The considered estimators are:
\begin{longlist}
\item[(1)] the Lasso
estimator with regularization parameter obtained by ten-fold
cross-validation;
\item[(2)] the  \textsc{mc}$+$ estimator of \citet{Zha10} with regularization parameter obtained by ten-fold
cross-valida\-tion;
\item[(3)] the  \textsc{scad} estimator  of \citet{FanLi01} with regularization parameter obtained by ten-fold
cross-validation.
\end{longlist}
The Lasso estimator is calculated using the \texttt{glmnet} package in R (\citeauthor{FriHasTib10}\break (\citeyear{FriHasTib10})).
The cross-validated \textsc{mc}$+$ and \textsc{scad} estimators are implemented in the \texttt{ncvreg} package in R (\cite{BreHua11}).

The performance of each of the four estimators, generically denoted
by $\hat \theta$ is measured by its prediction\vspace*{1pt} error $|\bX(\hat
\theta -\theta^*)|_2^2/n=\|{\sf f}_{\hat \theta} -{\sf
f}_{\theta^*}\|^2$.
 Moreover,
even though the estimation error $|\hat \theta -\theta^*|_2^2$ is
not studied above, we also report its values for a better comparison with other
simulation studies.

\subsubsection{Coordinatewise sparsity}
The vector $\theta^*$
is\break given by $\theta^*_j=\1(j \le S)$ for some fixed $S$ so that\vspace*{1pt}
$|\theta^*|_0=S$. Here, $\1(\cdot)$ denotes the indicator function.

\begin{figure*}[b]

\includegraphics{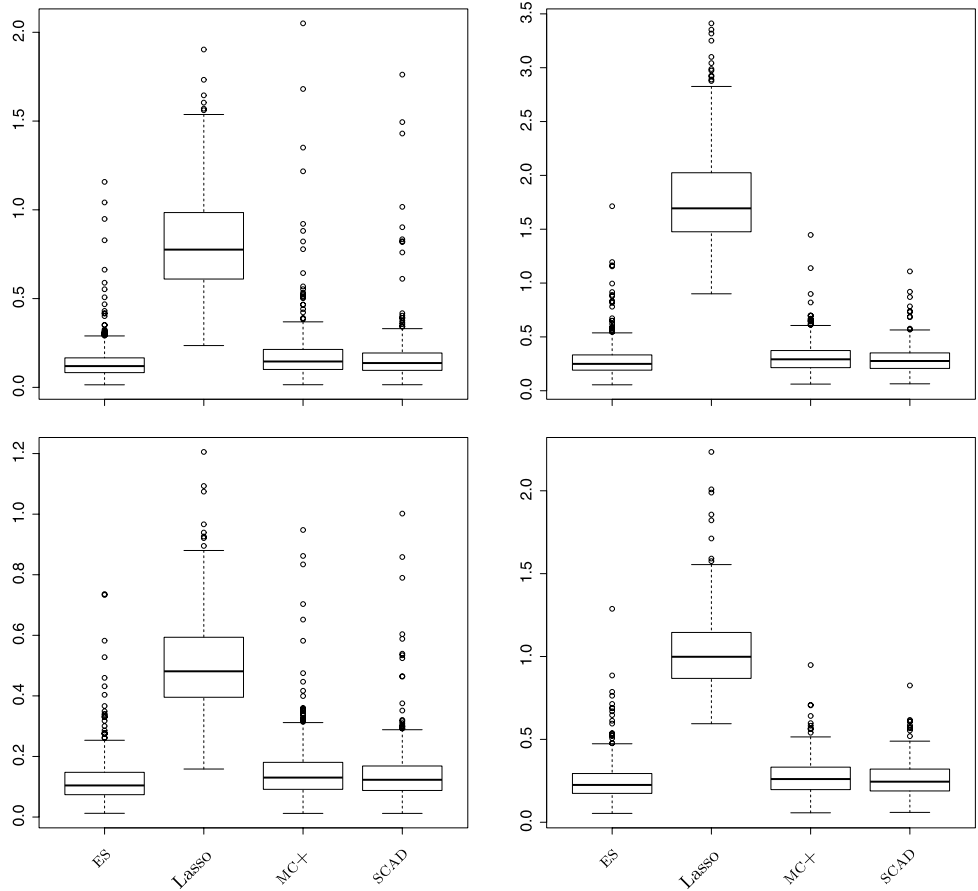}\vspace*{-3pt}

\caption{Boxplots of performance measure over 500 realizations for the \textsc{es}, Lasso, \textsc{mc}$+$ and \textsc{scad} estimators. \textit{Top:} estimation performance $|\hat \theta
-\theta^*|_2^2$. \textit{Bottom:} Prediction performance: $|\bX(\hat \theta
-\theta^*)|_2^2/n$. \textit{Left:} $(M,n,S)=(200,100,10)$. \textit{Right:} $(M,n,S)=(500,200,20)$.}
\label{FIG:boxplots}
\end{figure*}

We considered the cases $(n,M,S)\in \{(100, 200,
10),\allowbreak (200, 500, 20)\}$. The Metropolis approximation
$\tilde{\theta}^{\exp}_T$ was computed with $T_0=3000,  T=7000$,
which should be in the asymptotic regime of the Markov chain since
Figure~\ref{FIG:iter} shows that, on a typical example, the right
sparsity pattern is recovered after about 2000 iterations.

Figure~\ref{FIG:boxplots} displays comparative boxplots, and
Table~\ref{TAB:tabest} reports averages and standard deviations over
the 500 repetitions. In particular, it shows that \textsc{es}
outperforms the Lasso estimator and has performance similar to \textsc{mc}$+$ and \textsc{scad}.

Figure~\ref{FIG:iter} illustrates a typical behavior of the \textsc{es}
estimator for one particular realization of $\bX$ and $\xi$. For
better visibility,\vadjust{\goodbreak} both displays represent only the 50 first
coordinates of  $\tilde{\theta}^{\exp}_T$, with $T_0=3000,
T=7000$. The left-hand side display shows that the sparsity pattern
is well recovered and the estimated values are close to one. The
right-hand side display illustrates the evolution of the
intermediate parameter $\hat \theta_{\pP_t}$ for $t=1, \ldots,
5000$. It is clear that the Markov chain that runs on the
$M$-hypercube graph gets\break ``trapped'' in the vertex that corresponds
to the sparsity pattern of $\theta^*$ after only $2000$ iterations.
As a result, while the \textsc{es} estimator is not sparse itself, the
MH approximation to the \textsc{es} estimator may output a sparse
solution.

\begin{table}
\caption{Means and standard deviations of performance measures over 500 realizations for the \textsc{es}, Lasso, \textsc{mc}$+$ and \textsc{scad} estimators. \textit{Top:} estimation performance $|\hat \theta
-\theta^*|_2^2$. \textit{Bottom:} Prediction performance: $|\bX(\hat \theta
-\theta^*)|_2^2/n$}\label{TAB:tabest}
\begin{tabular*}{\columnwidth}{@{\extracolsep{\fill}}lcccc@{}}
\hline
$\bolds{(M,n,S)}$ & \textbf{\textsc{es}} & \textbf{Lasso} & \textbf{\textsc{mc}$+$} & \textbf{\textsc{scad}}\\
\hline
$(200,100, 10)$ &  0.14 & 0.82 & 0.18 & 0.17  \\
& (0.11) &    (0.28)  &   (0.17)   &  (0.15)  \\
[3pt]
$(500,200, 20)$ & 0.29 & 1.78 & 0.31 & 0.29 \\
&  (0.16) &     (0.43)  &  (0.14)  &  (0.12) \\
[6pt]
$(200,100, 10)$& 0.12 & 0.50 & 0.15 & 0.14   \\
& (0.08) &  (0.15)  & (0.10) &  (0.10)  \\
[3pt]
$(500,200, 20)$ & 0.25 & 1.02 & 0.27 & 0.26\\
& (0.11) &  (0.22)  & (0.11) & (0.10) \\
\hline
\end{tabular*}
\end{table}

\begin{Fusedsparsity*}
The vector $\theta^*$ is chosen piecewise
constant as follows. Fix an integer $S\ge 1$ such that $10S\le M$
and consider the blocks $I_1, \ldots, I_S$ defined by
\[
I_j=\{10(j-1)+1, \ldots, 10j\}, \quad  j=1,\dots,S.
\]
The vector $\theta^*$ is defined to take value $(-1)^j$ on $I_j,
j=1, \ldots, S$ and $1/2$ elsewhere.  We considered the cases
$(M, n,S)\in \{(200, 100, 10), (500, 200, 20)\}$ that are illustrated
in Figure~\ref{FIG:fused_est}. Note that in both cases, the vector $\theta^*$ is \textit{not} sparse.

The fused versions of Lasso, \textsc{mc}$+$ and \textsc{scad} are not
readily available in R, and we implement them as follows. Recall that
$D$ is the $M\times M$ matrix defined in Section~\ref{sub:priors}
by $(D\theta)_1=\theta_1$ and $(D\theta)_j=\theta_j-\theta_{j-1}$
for $j=2, \ldots, M$. The inverse $D^{-1}$ is the $M \times M$ lower
triangular matrix with ones on the diagonal and in the lower
triangle. To obtain the fused versions of Lasso, \textsc{mc}$+$ and \textsc{scad}, we simply run these algorithms on the design matrix $\bX
D^{-1}$ to obtain a solution $\hat \theta$. We then return the
vector $D^{-1}\hat \theta$ as a solution to the fused
problem.%

We report the boxplots of the two performance measures $|\bX(\hat
\theta -\theta^*)|_2^2/n$ and $|\hat \theta -\theta^*|_2^2$ in
Figure~\ref{FIG:boxplotsf}. It is clear that, in this example,
Exponential Screening outperforms the three other estimators.
Moreover, \textsc{mc}$+$ and \textsc{scad} perform particularly poorly in
the case $(M,n,S)=(500,200,20)$. Their output on a typical example
is illustrated in Figure~\ref{FIG:fused_est}. We can see that they
yield an estimator that takes only two values, thus missing most of
the structure of the problem. It seems that this behavior can be
explained by the fact that the estimators are trapped in a local
minimum close to zero.
\end{Fusedsparsity*}

\begin{appendix}
\section*{Appendix}\label{app}

\renewcommand{\theequation}{8.\arabic{equation}}
\setcounter{equation}{0}
The proof of~\eqref{EQ:TH:ERM:UB} is standard, and similar results
have been formulated in the literature for various other setups. We
give it here for the sake of completeness. From the definition of
the empirical risk minimizer~$\hat f^\textsc{erm}$, we have
\[
\hat R_n(\hat f^\textsc{erm})\le \hat R_n(f^*),
\]
where $f^*$ is any minimizer of the true risk $R(\cdot)$ over~$\cH$.
Simple algebra yields
\[
R(\hat f^\textsc{erm})\le R(f^*)+2\E\langle \hat f^\textsc{erm}
-f^*, \bY - \eta \rangle,
\]
where for two functions $f,g$ from $\cX$ to $\R$ we set $ \langle f,
g\rangle =\frac{1}{n}\sum_{i=1}^nf(x_i)g(x_i)$.  Next, observe
that
\begin{eqnarray*}
\E\langle \hat f^\textsc{erm} -f^*, \bY - \eta \rangle&\le& \E\max_{f
\in \cH} \langle f-f^*, \bY - \eta \rangle\\
& \le& 2\sigma
\sqrt{\frac{2\log M}{n}},
\end{eqnarray*}
where we used the fact that $\|f^*-f\| \le 2$ for any \mbox{$f\in \cH$}, and the inequality
$\E[\max_{1\le i\le M}a_i^\top \xi]\le \sigma\cdot\allowbreak \sqrt{2\log M}$ valid for any $a_1, \ldots, a_n \in \R^n$, $|a_i|_2 \le 1$, where $\xi=(\xi_1, \ldots, \xi_n)^\top$.

\renewcommand{\thefigure}{\arabic{figure}}
\begin{figure*}

\includegraphics{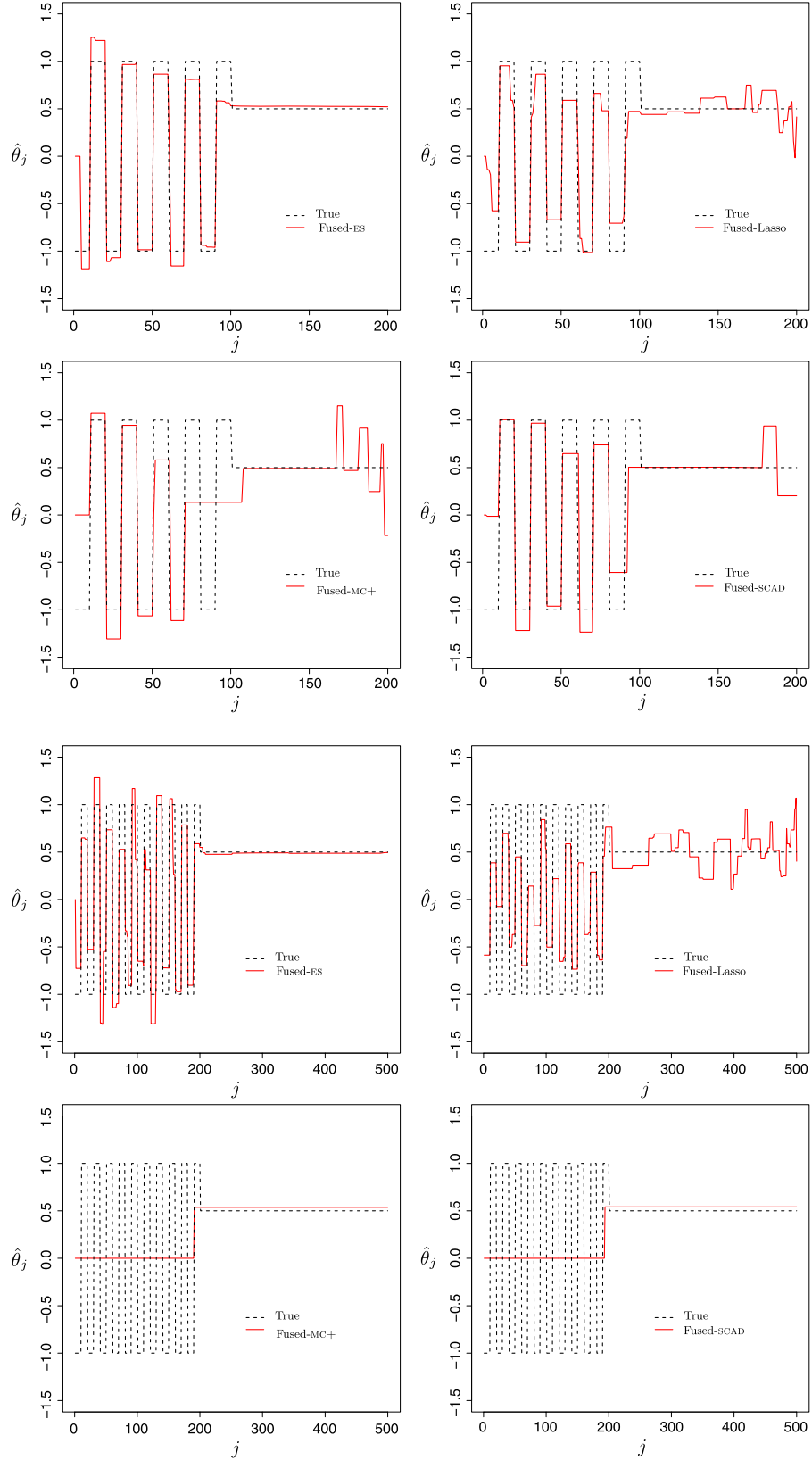}

\caption{Typical realizations of the fused estimators in the cases $(M,n,S)=(200,100,10)$ (top) and $(M,n,S)=(500,200,20)$ (bottom).} \label{FIG:fused_est}
\end{figure*}

\begin{figure*}

\includegraphics{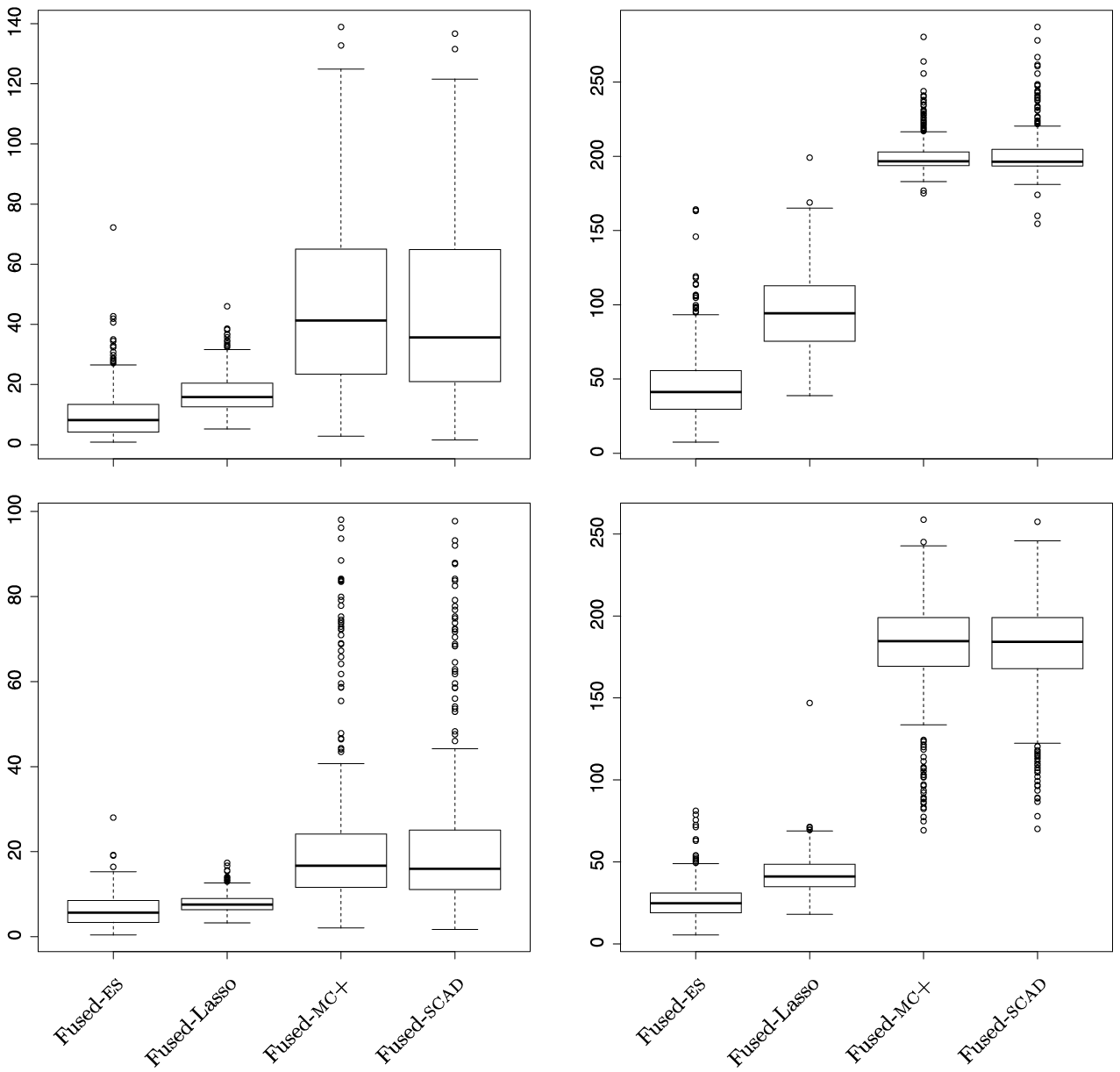}

\caption{Boxplots of performance measure over 500 realizations for the Fused-\textsc{es}, Fused-Lasso, Fused-\textsc{mc}$+$ and Fused-\textsc{scad} estimators. \textit{Top:} estimation performance $|\hat \theta
-\theta^*|_2^2$. \textit{Bottom:} Prediction performance: $|\bX(\hat \theta
-\theta^*)|_2^2/n$. \textit{Left:} $(M,n,S)=(200,100,10)$. \textit{Right:} $(M,n,S)=(500,200,20)$.}  \label{FIG:boxplotsf}
\end{figure*}

We now turn to the proof of~\eqref{EQ:TH:ERM:LB}. Consider the
random matrix $\mathbb{X}$ of size $n\times M$ such that its elements
$\mathbb{X}_{i,j}, i=1, \ldots, n,  j=1, \ldots, M$ are i.i.d. Rademacher
random variables, that is, random variables taking values $1$ and $-1$
with probability $1/2$.
%
%
Moreover, assume that
\begin{equation}
\label{EQ:cond_rip}
\frac{2}{n}\log \biggl(1+\frac{eM}{2} \biggr) <C_1,
\end{equation}
for some positive constant $C_1<1/2$. Note that~\eqref{EQ:cond_rip} follows
from~\eqref{EQ:C0} if $C_0$ is chosen small enough. Theorem~5.2 in
\citet{BarDavDeV08} (see also Section~5.2.1 in \cite{RigTsy11})
entails that if~\eqref{EQ:cond_rip} holds for $C_1$ small enough,
then there exists a nonempty set $\cM$ of matrices obtained as
realizations of the matrix $\X$ that enjoy the following weak
restricted isometry (\textsc{wri}) property. For any $X \in \cM$, there
exists constants $\underline{\kappa}\ge \bar \kappa>0$, such that
for any $\lambda \in \R^M$ with at most $2$ nonzero coordinates,
\begin{equation}
\label{EQ:rip}
\underline{\kappa}^2|\lambda|_2^2 \le \frac{|X \lambda|_2^2}{n} \le \bar \kappa^2 |\lambda|_2^2,
\end{equation}
when~\eqref{EQ:cond_rip} is satisfied. For $X \in \cM$, let $\phi_1,
\ldots, \phi_M$ be any functions on $\cX$ satisfying
\[
\phi_j(x_i)=X_{i,j},\quad   i=1, \ldots, n, j=1, \ldots, M,
\]
where $X_{i,j}$ are the entries of $X$. Note that $\|\phi_j\|=1$
since $X_{i,j}\in \{-1,1\}$.

Fix $\tau>0$ to be chosen later, and set
\[
f_j= \tau (1 +\alpha )\phi_j , \quad  j=1,\dots,M,
\]
where we set for brevity $\alpha=(\sigma/3) \sqrt{\frac{\log M}{\bar
\kappa^2n}}$. Moreover, consider the functions
\[
\eta_j = \tau\alpha\phi_j ,\quad   j=1,\dots,M.
\]
Using~\eqref{EQ:C0} we choose $\tau$ small enough to ensure that $\|\eta_j\|\le 1$ and $\|f_j\|\le 1$ for any $j=1, \ldots, M$.

We write $R_j(\cdot)$ to denote the risk function $R(\cdot)$
when $\eta=\eta_j$ in~\eqref{EQ:model}. It is easy to check that
\begin{equation}
\label{EQ:prTH1_1}
\min_{f \in \cH}R_j(f)=R_j(f_j)=\|f_j -\eta_j\|^2 .
\end{equation}
As it is customary in the proof of minimax lower bounds, we reduce
our estimation problem to a testing problem as follows. Let $\psi
\in \{1, \ldots, M\}$ be the random variable, or \textit{test},
defined by $\psi=j$ if and only if $\hat S_n=f_j$. Then, $\psi\neq
j$ implies that there exists $k \neq j$ such that $\hat S_n=f_k$, so
that
\begin{eqnarray*}
&&\|\hat S_n -\eta_j\|^2-\|f_j -\eta_j\|^2 \\
&&\quad=\|f_k -f_j\|^2   + 2
\langle f_k -f_j, f_j -\eta_j\rangle \\
&&\quad= \tau^2(1+\alpha)^2\|\phi_j-\phi_k\|^2\\
&&\qquad{} +
2\tau^2(1+\alpha)(\langle\phi_j,\phi_k\rangle-1)\\
&&\quad\ge \tau^2 \alpha \|\phi_j-\phi_k\|^2.
\end{eqnarray*}
From~\eqref{EQ:rip}, we find that $\|\phi_j -\phi_k\|^2\ge 2 \underline{\kappa}^2$ so that
\[
\|\hat S_n -\eta_j\|^2-\|f_j -\eta_j\|^2 \ge
\frac{2\tau^2\underline{\kappa}^2\sigma}{3\bar \kappa}
\sqrt{\frac{\log M}{n}}\triangleq \nu_{n,M}.
\]
Therefore, we conclude that $\psi \neq j$ implies that
\[
R_j(\hat S_n)-\min_{f \in \cH}R_j(f) \ge \nu_{n,M}.
\]
Hence,
\begin{eqnarray}
\label{EQ:minimax_pr}
&&\max_{1\le j\le M}P_j \Bigl\{R_j(\hat S_n)-
\min_{f \in \cH}R_j(f) \ge \nu_{n,M} \Bigr\}\nonumber\\[-8pt]\\[-8pt]
&&\quad \ge \inf_{\psi}
\max_{1\le j\le M}P_j(\psi \neq j) ,\nonumber
\end{eqnarray}
where the infimum is taken over all tests taking values in $\{1,
\ldots, M\}$, and $P_j$ denotes the joint distribution of $Y_1,
\ldots, Y_n$ that are independent Gaussian random variables with
mean $\eta_j(x_i)$, respectively. It follows
from \citeauthor{Tsy09} [(\citeyear{Tsy09}), Proposition~2.3 and Theorem~2.5] that if for any
$1\le j,k\le M$, the Kull\-back--Leibler divergence between $P_j$ and
$P_k$ satisfies
\begin{equation}
\label{EQ:KL}
\cK(P_j, P_k)<\frac{\log M}{8} ,
\end{equation}
then there exists a constant $C>0$ such that
\begin{equation}
\label{EQ:inf_psi}
\inf_{\psi} \max_{1\le j\le M}P_j(\psi \neq j)\ge C .
\end{equation}
To check~\eqref{EQ:KL}, observe that, choosing $\tau\le 1$ and
applying~\eqref{EQ:rip}, we get
\begin{eqnarray*}
\cK(P_j, P_k)&=&\frac{n}{2\sigma^2}\|\eta_j-\eta_k\|^2=
\frac{\tau^2\log M}{18\bar \kappa^2}\|\phi_j -\phi_k\|^2\\
&<&\frac{\log
M}{8} .
\end{eqnarray*}
Therefore, in view of~\eqref{EQ:minimax_pr} and~\eqref{EQ:inf_psi},
we find, using the Markov inequality, that for any selector $\hat
S_n$,
\begin{eqnarray*}
\max_{1\le j\le M}E_j \Bigl[R_j(\hat S_n)- \min_{f \in
\cH}R_j(f) \Bigr] &\ge& C\nu_{n,M}\\
& =&C_*\sigma \sqrt{\frac{\log
M}{n}} ,
\end{eqnarray*}
where $E_j$ denotes the expectation with respect to~$P_j$.
\end{appendix}

\section*{Acknowledgments}
The first author is supported in part by the NSF DMS-09-06424, DMS-10-53987.


\end{document}